\documentclass{amsart}
\usepackage[utf8]{inputenc}
\usepackage{enumerate}
\usepackage{enumitem}
\usepackage{graphicx}
\usepackage{wrapfig}
\usepackage{subcaption}
\usepackage{listings}
\usepackage{color}
\definecolor{mygreen}{RGB}{28,112,30} 
\definecolor{mylilas}{RGB}{170,55,241}
\definecolor{myblue}{RGB}{20,30,171}
\definecolor{myred}{RGB}{200,10,30}
\usepackage{microtype}
\usepackage{multirow}
\usepackage{array}
\newcounter{rowcount}
\usepackage[foot]{amsaddr}
\usepackage{amsmath}
\usepackage{amsfonts}
\usepackage{amssymb}
\usepackage{amscd}
\usepackage{amstext}
\usepackage{amsbsy}
\usepackage{amsopn}
\usepackage{amsthm}
\usepackage{thmtools}
\usepackage{amsxtra}
\usepackage{float}
\usepackage{bbm}
\usepackage{url}
\usepackage{hyperref}
\hypersetup{
colorlinks=true,
linkcolor=myred,
citecolor=blue}

\newtheoremstyle{mytheoremstyle} 
  {12pt}                    
  {\topsep}                    
  {\slshape}                   
  {}                           
  {\bfseries}                   
  {.}                          
  {.5em}                       
  {}  

\newtheoremstyle{mydefstyle} 
  {12pt}                    
  {12pt}                    
  {}                   
  {}                           
  {\bfseries}                   
  {.}                          
  {.5em}                       
  {}  

  \newtheoremstyle{plainsl}%
  	{}
  	{\topsep}
  	{\slshape} 
  	{}
  	{\normalfont\bfseries}
  	{.}
  	{ }
  	{}

\theoremstyle{plainsl}
\newtheorem{thm}{Theorem}[section]
\newtheorem*{thm*}{Theorem}

\newtheorem{cor}[thm]{Corollary}
\newtheorem{lem}[thm]{Lemma}
\newtheorem*{lem*}{Lemma}

\theoremstyle{remark}

\theoremstyle{mydefstyle}

\newtheorem{xmpl}[thm]{Example}
\newtheorem*{note*}{Remark}

\newcommand{\Z}{\mathbb{Z}}
\newcommand{\N}{\mathbb{N}}
\newcommand{\R}{\mathbb{R}}

\newcommand{\Cx}{\mathbb{C}}
\newcommand{\A}{\mathcal{A}}
\newcommand{\M}{\mathcal{M}}
\newcommand{\1}{\mathbbm{1}}

\newcommand{\C}{\mathcal{C}}
\renewcommand{\S}{\mathcal{S}}

\renewcommand{\phi}{\varphi}
\newcommand{\seq}{\subseteq}

\newcommand{\B}{\mathcal{B}}
\newcommand{\D}{\mathcal{D}}
\newcommand{\sm}{\backslash}

\newcommand\diff{\mathbin{\mkern-1.5mu\setminus\mkern-1.5mu}}

\renewcommand\qed{%
	\ifmmode\eqno\sqr53
	\else\nolinebreak\ \hfill\sqr53\medbreak\fi}
\renewcommand\proof{\noindent\textsl{Proof. }}
\newcommand\sqr[2]{{\vbox{\hrule height.#2pt
    \hbox{\vrule width.#2pt height#1pt \kern#1pt
        \vrule width.#2pt}\hrule height.#2pt}}}

\def\pmat#1{{\begin{pmatrix}#1\end{pmatrix}}}

\title{On state transfer in Cayley graphs for abelian groups}
\author{Arnbjörg Soffía Árnadóttir$^{1^*}$}
\address{$^{1^*}$Corresponding author. Department of Applied Mathematics and Computer Science, Technical University of Denmark, DK-2800
Lyngby, Denmark.  e-mail:  {\tt sofar@dtu.dk}}

\author{Chris Godsil$^2$}
\address{$^2$Department of Combinatorics and Optimization, University of Waterloo,
200 University Avenue West, Waterloo, ON, Canada N2L 3G1.  e-mail:  {\tt cgodsil@uwaterloo.ca}}

\thanks{Both authors acknowledge the support of C.\ Godsil's NSERC (Canada), Grant No.\ RGPIN-9439}

\begin{document}
\renewcommand{\itshape}{\slshape}

\begin{abstract}
  In this paper, we characterize perfect state transfer in Cayley graphs for abelian groups that have a cyclic Sylow-2-subgroup. This generalizes a result of Ba\v si\'c from 2013 where he provides a similar characterization for Cayley graphs of cyclic groups.
\end{abstract}

\maketitle
\setlength{\parskip}{0cm}
\tableofcontents
\setlength{\parskip}{0.2cm}

\section{Introduction}

State transfer on graphs has been of interest to physicists and mathematicians for the past twenty years, starting with papers by Bose in 2003 \cite{bose2003} and Christandl et al in 2004 \cite{christandl2004}.
Given a graph $X$ with adjacency matrix $A$, the \textit{continuous-time quantum walk} on $X$ is given by the matrices $U_A(t)=e^{itA}$ for $t\in \R$. The matrix $U_A(t)$ is called the \textit{transition matrix} of the quantum walk at time $t$. We say that there is \textit{perfect state transfer} from vertex $u$ to vertex $v$ at time $t$ if $|U_A(t)_{u,v}|=1$.

Continuous-time quantum walks were first introduced by Farhi and Gutmann \cite{farhi-gutmann} in the context of search on graphs. Quantum walks had then already been used for graph search in discrete time, for instance by Grover \cite{grover}. It turns out that both in discrete and continuous time, quantum walks perform significantly better than classical random walks. In 2003, Childs et al \cite{childs2003} constructed a search problem that can be solved efficiently with a continuous-time quantum walk and proved that no classical algorithm can solve this problem in subexponential time.

In the same year, Bose used the formulation by Farhi and Gutmann to study information transmission in quantum spin chains, and this is where the notion of state transfer arises. For more on the applications of perfect state transfer in physics, we refer the reader to a review on the subject by Kay \cite{kay2010}.

In quantum physics, the No-Cloning Theorem states that it is impossible to create a new copy of a quantum state. This is why perfect state transfer is useful: it is a way to effectively transfer a quantum state when it cannot be copied.
Unfortunately however, perfect state transfer in graphs is quite rare. In \cite{godsil2012when}, Godsil shows that for a fixed $k$, there are only finitely many connected graphs with maximum degree $k$ that admit perfect state transfer. Thus, many results about perfect state transfer show non-existence in certain graphs rather than providing examples.

Perfect state transfer between vertices in a graph implies that these vertices are ``similar'' in some sense, thus it seems natural to look for this behaviour in vertex-transitive graphs. In particular, it has been popular to look at Cayley graphs.
In 2008, Bernasconi, Godsil and Severini considered state transfer in cubelike graphs (Cayley graphs of elementary abelian 2-groups) \cite{bernasconi}. Their results imply that most cubelike graphs have perfect state transfer. More results on state transfer in cubelike graphs were later proved by Cheung and Godsil in 2010 \cite{cheung} and Chan in 2013 \cite{chan2013complex}.

Another important class of Cayley graphs that is of interest in the context of state transfer is that of circulants (Cayley graphs for cyclic groups). This was first considered by Saxena, Severini and Shparlinski in 2007 \cite{saxena2007}, where they proved some necessary conditions. In the years 2009--2013, a series of four papers was written by various subsets of $\{\text{Ba\v si\'c, Petkovi\'c, Stevanovi\'c}\}$, on perfect state transfer in circulants \cite{basic2009first, basic2009some, basic2010nonsq, basic2013char}. In the last one of these, Ba\v si\'c completely characterizes the connection sets of circulants having perfect state transfer.

In this paper, we generalize Ba\v si\'c's characterization to Cayley graphs of abelian groups having a cyclic Sylow-$2$-subgroup. We will call such graphs $2$-circulants. While Ba\v si\'c et al used number theory in their proofs, we approach the topic from a group theoretical perspective, thus not only providing a generalization, but also completely different proofs of Ba\v si\'c's results.
Our main result is the following.
\begin{thm*}[Theorem \ref{thm:final}]
  Let $G$ be an abelian group of order $2^dm$ where $m$ is odd and suppose it has a cyclic Sylow-2-subgroup. Let $X=X(G,\C)$ be a Cayley graph. If $d=0$, there is no perfect state transfer on $X$ and if $d=1$, there is perfect state transfer if and only if $X$ is a matching. Suppose $d\geq 2$, let $a$ be the unique element of order two, and $b,-b$ the unique pair of elements of order four. Denote by $\C_k$ the set of elements in $\C$ with order $2^km'$ where $m'$ is odd. Then $X$ has perfect state transfer if and only if
  \begin{enumerate}[label = (\alph*)]
    \item $\C$ is power-closed,
    \item either $a$ or $b$ is in $\C$ but not both,
    \item $\C_0=4(\C_2\diff\{-b,b\})$, and
    \item $\C_1\diff\{a\} = 2(\C_2\diff\{-b,b\})$.
  \end{enumerate}
\end{thm*}
We also describe the spectrum of a $2$-circulant having perfect state transfer (which again is consistent with Ba\v si\'c's description of the spectrum of circulants) and give the time at which perfect state transfer must occur in such graphs. In the last section, we propose that further generalization might be possible, but this is currently work in progress.

\section{Preliminaries}\label{sec:prelim}

Let $A$ be a hermitian matrix. The \textit{continuous-time quantum walk} on $A$ at time $t\in\R$ is given by the matrix
\[U_A(t):=e^{itA} = \sum_{n\geq0}\frac{(it)^n}{n!}A^n.\]
We call $U_A(t)$ the \textit{transition matrix} of $A$ at time $t$, it is a unitary matrix for all $t$. We will generally take $A$ to be the adjacency matrix of a graph, $X$. Most of our graphs are simple and undirected, but we will need the notions of weighted graphs and graphs with loops later on. In all cases, the adjacency matrix is real and symmetric, and thus hermitian. In this case we talk about a quantum walk on the graph $X$, and write $U_X(t)$, or even $U(t)$ if the graph is clear from the context.

Let $X$ be a simple, undirected graph with adjacency matrix $A(X)=A$ and let $u$ and $v$ be vertices of $X$. We say that there is \textit{perfect state transfer} from $u$ to $v$ at time $\tau$ with phase factor $\lambda$ if $U(\tau)\mathbf{e}_v=\lambda \mathbf{e}_u$, where $\mathbf{e}_x$ denotes the standard basis vector indexed by $x$. Equivalently, $|U(\tau)_{u,v}|=1$.
We say that the vertex $u$ is \textit{periodic} at time $\tau$ with phase factor $\lambda$ if $U(\tau)\mathbf{e}_u=\lambda \mathbf{e}_u$, equivalently, $|U(\tau)_{u,u}|=1$. If every vertex is periodic at time $\tau$, we call the graph $X$ \textit{periodic} at time $\tau$.
The minimum time at which a vertex or a graph is periodic is called the \textit{period} of the vertex or graph, respectively.

We refer to the eigenvalues and eigenvectors of the matrix $A$ as the eigenvalues and eigenvectors of the graph $X$. Let $\theta_0,\dots,\theta_n$ be the (distinct) eigenvalues of $X$ and let $E_r$ denote the orthogonal projection onto the eigenspace of $\theta_r$. Then each $E_r$ 
satisfies $E_r^2=E_r$ and if $s\neq r$, we have $E_rE_s=0$. Furthermore, $\sum_{r=0}^n E_r=I$ and
\[A = \sum_{r=0}^n \theta_rE_r.\]
This is called the \textit{spectral decomposition} of $A$ and the matrices $E_r$ are its \textit{spectral idempotents.} An important and useful fact is that if $f$ is a univariate function defined on the spectrum of $A$, then
\[f(A) = \sum_{r=0}^nf(\theta_r) E_r.\]
In particular, if $U(t)$ is the transition matrix of a quantum walk on $X$, then
\[U(t) = \sum_{r=0}^n e^{it\theta_r}E_r.\]
Observe that an immediate consequence of this is that if a graph has integer eigenvalues, then it is periodic at time $\tau=2\pi$.

Let $G$ be a group and let $\C\subseteq G$ be an inverse-closed subset. We define the \textit{Cayley graph} for $G$ with respect to $\C$ to be the graph with vertex set $G$, in which vertices $g$ and $h$ are adjacent if and only if $hg^{-1}\in\C$. We denote this graph by $X(G,\C)$ and refer to $\C$ as its \textit{connection set}. Note that the graph has a loop on every vertex if and only if $e\in \C$ (otherwise it has no loops).
We will use Cayley graphs with loops in some of our proofs, however, it is easy to verify that there is perfect state transfer on the matrix $A$ if and only if there is perfect state transfer on $A+I,$ thus the loops do not play any role in our final results.

A Cayley graph for an abelian group is called a \textit{translation graph} and a Cayley graph for a cyclic group is called a \textit{circulant}. In this paper we will focus on Cayley graphs for abelian groups that have a cyclic Sylow-2-subgroup. Such graphs will be called \textit{$2$-circulants}.

We define a \textit{weighted Cayley graph} as a Cayley graph $X(G,\C)$ together with a function, $\omega:\C\to\Z$ such that the fibre of each element in $\Z$ is inverse-closed. The \textit{weighted adjacency matrix} of this graph is given by
\[\sum_{n\in \omega(\C)}nA(X(G,\omega^{-1}(n))).\]
In particular, a \textit{signed Cayley graph} is a weighted Cayley graph for which the image of $\omega$ is $\{\pm 1\}$. Denote by $\C^+$ and $\C^-$ the subset of elements of $\C$ with positive and negative sign, respectively. Then the \textit{signed adjacency matrix} of $X(G,\C)$ is given by $A(X(G,\C^+))-A(X(G,\C^-))$.

Given graphs $X$ and $Y$ with adjacency matrices $A(X)$ and $A(Y)$, we define their \textit{Cartesian product} $X\square Y$, as the graph that has adjacency matrix $A(X)\otimes I+ I\otimes A(Y)$, where $\otimes$ denotes the Kronecker product of matrices. Equivalently, this is the graph with vertex set $V(X)\times V(Y)$, in which $(x_1,y_1)\sim(x_2,y_2)$ are adjacent if and only if either $x_1=x_2$ and $y_1\sim y_2$, or $x_1\sim x_2$ and $y_1=y_2$.
We define the \textit{direct product} of $X$ and $Y$, denoted by $X\times Y$, as the graph with adjacency matrix $A(X)\otimes A(Y)$. This is the graph with vertex set $V(X)\times V(Y)$ and $(x_1,y_1)\sim (x_2,y_2)$ if and only if both $x_1\sim x_2$ and $y_1\sim y_2$.

Clearly, the matrices $A(X)\otimes I$ and $I\otimes A(Y)$ commute, and now it is not too hard to see that $U_{X\square Y}(t) = U_X(t)\otimes U_Y(t)$.

\section{Perfect state transfer}\label{sec:pst}
In this section we will give some preliminary results about perfect state transfer and periodicity in graphs.

\begin{lem}[{\cite[Lemma 1.1]{godsil2012state}}]\label{lem:pstsym}
  Let $u$ and $v$ be vertices of a graph $X$. If there is perfect state transfer from $u$ to $v$ at time $\tau,$ then there is perfect state transfer from $v$ to $u$ at time $\tau$ and $u$ and $v$ are both periodic at time $2\tau.$\qed
\end{lem}
By Lemma \ref{lem:pstsym}, we can talk about perfect state transfer occurring between two vertices, rather than from one to another.

The following lemma on so-called monogamy of perfect state transfer was first observed by Kay in 2011 \cite[Section D]{kay2011}. A proof can also be found in \cite{godsil2012state}.

\begin{lem}[{\cite[Corollary 3.2]{godsil2012state}}]
  If there is perfect state transfer from $u$ to $v$ in $X$ and also from $u$ to $w$, then $v=w$. \qed
\end{lem}

In 2011, Godsil gave a characterization of periodicity of graphs in terms of their eigenvalues.

\begin{lem}[{\cite[Corollary 3.3]{godsil2011periodic}}]\label{lem:periodic-int}
  A graph $X$ is periodic if and only if either
  \begin{enumerate}
    \item the eigenvalues of $X$ are integers, or
    \item the eigenvalues of $X$ are rational multiples of $\sqrt\Delta$, for some fixed, square-free integer $\Delta.$\qed
  \end{enumerate}
\end{lem}

We now consider what happens in vertex-transitive graphs in particular. Note that all Cayley graphs are vertex-transitive: if $X=X(G,\C)$ is a Cayley graph, then $G$ acts regularly (in particular transitively) on $X$ as a group of automorphisms.

Let $X$ be a vertex-transitive graph. If there is perfect state transfer between two vertices in $X$ at time $\tau$, then every vertex of $X$ is involved in perfect state transfer at time $\tau$ and the graph is periodic at time $2\tau$. This is a result of the following theorem of Godsil.

\begin{thm}[{\cite[Theorem 6.1]{godsil2012state}}]\label{thm:pstvxtr}
  Let $X$ be a vertex-transitive graph and let $u$ and $v$ be vertices of $X$. If there is perfect state transfer from $u$ to $v$ at time $\tau$, then $U(\tau)$ is a scalar multiple of a permutation matrix with order two and no fixed points.\qed
\end{thm}

The theorem implies that a vertex-transitive graph admitting perfect state transfer has an even number of vertices. In particular, a Cayley graph for a group of odd order has no perfect state transfer. In fact, Theorem \ref{thm:pstvxtr} holds for integer-weighted vertex-transitive graphs, and so weighted Cayley graphs for groups of odd order have no perfect state transfer.

Theorem \ref{thm:pstvxtr} further implies that a vertex transitive graph with perfect state transfer is periodic, and now the next lemma follows from Lemma \ref{lem:periodic-int} and the fact that a vertex-transitive graph has an integer eigenvalue, namely its degree.
\begin{lem}\label{lem:integrality}
  Let $X$ be a vertex-transitive graph. If $X$ admits perfect state transfer, then all its eigenvalues are integers.\qed
\end{lem}
We call a graph with integer eigenvalues \textit{integral}.

\section{Association schemes}
In this section we will give a brief introduction to association schemes, particularly group schemes. We will make many claims without proofs here, but details can be found in \cite[Chapter 2]{brouwer1989} or \cite[Chapters 3 and 11]{godsil2016EKR}. We conclude the section with a bit of character theory, but for more on representations and characters of groups we refer the reader to \cite{james2001reps}.

Let $J$ denote the $n\times n$ all-ones matrix. An \textit{association scheme with $d$ classes} is a set of $n\times n$ matrices, $\A = \{A_0,\dots,A_d\}$ with entries in $\{0,1\}$ such that
\begin{enumerate}
  \item $A_0=I$ and $\sum_{r=0}^d A_r = J$,
  \item $A_r^T\in \A$ for all $r$,
  \item $A_rA_s=A_sA_r$ for all $r,s$,
  \item $A_rA_s$ lies in the span of $\A$ for all $r,s$.
\end{enumerate}

The span of $\A$ is a commutative algebra, $\Cx[\A]$, called the \textit{Bose-Mesner algebra} of the association scheme and any $\{0,1\}$-matrix in this algebra is called a \textit{Schur idempotent} of $\Cx[\A]$. The elements in the scheme, $A_0,\dots, A_d$, are the \textit{minimal Schur idempotents} and every Schur idempotent is a sum of some minimal Schur idempotents. If $\B$ is an association scheme such that every element of $\B$ is a Schur idempotent of $\Cx[\A]$, we say that $\B$ is a \textit{subscheme} of $\A$.

The minimal Schur idempotents of $\A$ form a basis for the Bose-Mesner algebra. It can be shown that the algebra has another basis, $E_0,\dots,E_d$ of matrix idempotents that sum to the identity matrix. Further, there are scalars, $p_r(s)$ and $q_s(r)$, for $r,s = 0,\dots,d$ such that
\[A_r = \sum_{s=0}^d p_r(s)E_s\quad\text{and}\quad E_s = \frac1n\sum_{r=0}^d q_s(r)A_r.\]
Define matrices, $P=(p_r(s))_{s,r}$ and $Q=(q_s(r))_{r,s}$. We call $P$ the \textit{matrix of eigenvalues} of the scheme $\A$ (the scalars $p_r(s)$ are eigenvalues of the $A_r$) and $Q$ the \textit{matrix of dual eigenvalues} of $\A$, and we have $PQ=nI.$ We say that an association scheme is \textit{formally self-dual} if $Q = \overline{P}.$


Suppose $G$ is a group of order $n$, with conjugacy classes $C_0,\dots,C_d$, where $C_0=\{e\}$. Define $A_r$ to be the $n\times n$ matrix indexed by the group elements where the $(g,h)$-entry is one if $hg^{-1}\in C_r$ and zero otherwise. The set $\A:=\{A_0,\dots, A_d\}$ is called the \textit{conjugacy class scheme} on $G$ and it is an association scheme. A detailed discussion on conjugacy class schemes can be found in \cite[Chapter 11]{godsil2016EKR}.
Observe that $A_r$ is the adjacency matrix for the (possibly directed) Cayley graph $X(G,\C_r)$, and in fact, every symmetric Schur idempotent in $\Cx[\A]$ is the adjacency matrix of some Cayley graph of $G$. We refer to these graphs as the graphs in the scheme.

Now let $G$ be an abelian group. Then its conjugacy classes all have size one, and we refer to the conjugacy class scheme as the \textit{abelian group scheme} of $G$. A subscheme of an abelian group scheme is called a \textit{translation scheme,} and now we have an alternative definition of a translation graph: it is a graph in a translation scheme. Let $\A$ be the abelian group scheme of $G$. We are interested in one particular subscheme of $\A$. Consider the subspace of the Bose-Mesner algebra, $\Cx[\A]$, spanned by the matrices with rational entries and rational eigenvalues. It turns out that this is the Bose-Mesner algebra of a subscheme, $\B$ of $\A$ and the integer matrices in $\Cx[\B]$ have only integer eigenvalues. We call $\B$ the integral translation scheme of $G$. Every integral translation graph of $G$ lies in the integral translation scheme of $G$.

Define the \textit{characters} of the abelian group $G$ to be the group homomorphisms from $G$ to the multiplicative group of the complex numbers (note that this definition coincides with the definition of an irreducible character of a representation of the group.) The characters form a group under pointwise multiplication, called the \textit{character group} of $G$ and denoted by $G^*$. This group is isomorphic to $G$.

Let $\chi$ be a character of $G$, and let $X=X(G,\C)$ be a translation graph with adjacency matrix $A$. Define
\[\chi(\C):= \sum_{g\in\C}\chi(g).\]
Then using the fact that $\chi$ is a homomorphism, it is easy to show that $A\chi = \chi(\C)\chi$, and so $\chi$ is an eigenvector of $X$ with eigenvalue $\chi(\C).$ Therefore, all translation graphs of $G$ have the same eigenvectors and we see that if $\C$ and $\D$ are disjoint subsets of $G$ such that the Cayley graphs $X(G,\C)$ and $X(G,\D)$ have eigenvalues $\theta$ and $\lambda$, respectively, both for the eigenvector $\chi$, then $\theta+\lambda$ is an eigenvalue of $X(G,\C\cup\D)$.

If $X$ is an integral translation graph, its adjacency matrix is a Schur idempotent of the integral translation scheme, $\A=\{A_0,\dots,A_d\}$ and thus a sum of some minimal Schur idempotents. Let $P$ be the matrix of eigenvalues of $\A$ and let $x$ be the vector of length $d+1$ with $1$ in the $r$-th entry if $A_r$ is in this sum, and zero otherwise. Then the eigenvalues of $X$ are the entries of $Px$.

Define an equivalence relation on $G$ by letting $g\approx h$ for $g,h\in G$ if and only if $\langle g\rangle = \langle h\rangle$, that is if $g$ and $h$ generate the same cyclic subgroup of $G$. We then say that $g$ and $h$ are \textit{power-equivalent} and refer to the equivalence classes of this relation as \textit{power classes}. Denote the power class of $g\in G$ by $[g]$. We call a subset of $G$ \textit{power-closed} if it is a union of power classes. The following lemma is due to Bridges \& Mena \cite{Bridges1982}.

\begin{lem}[{\cite[Corollary 2.5]{Bridges1982}}]
  The translation graph $X=X(G,\C)$ is integral if and only if $\C$ is power-closed.\qed
\end{lem}

This lemma can be used to show that the integral translation scheme of an abelian group is formally self-dual. A proof of this can be found in \cite[Theorem 4.7.3]{ArnadottirThesis}.


\section{Integral translation graphs}\label{sec:translation}
In this section, we continue our discussion of integral translation graphs. In what follows, $G$ will be an abelian group written additively.

We call $X(G,\C)$ an \textit{integral signed Cayley graph} if both $X(G,\C^+)$ and $X(G,\C^-)$ are integral. In this case, both the signed and unsigned adjacency matrices have integer eigenvalues.

\begin{lem}
  \label{lem:oddeval}
  If $G$ is an abelian group of odd order, then any non-empty, integral signed Cayley graph for $G$ has an odd eigenvalue.
\end{lem}
\proof
  Let $X$ be an integral signed Cayley graph for the abelian group $G$ of odd order $n$. Let $\A=\{A_0,\dots,A_d\}$ be the integral translation scheme for $G$ and let $P$ be the matrix of eigenvalues of $\A$. Since $\B$ is self-dual we have $P\overline{P}=nI$. Further, $P$ is real, so its determinant is real, and thus $\det(P)=\det(\overline{P}).$
  Then
  \[\det(P)^2 =\det(P)\det(\overline{P}) =\det(P\overline{P}) =\det(nI)=n^{d+1}.\]
  Therefore, $\det(P)$ is odd implying that $P$ is invertible modulo $2$.

  Now recall that $X$ lies in the scheme $\A$, so there is a $\{0,1\}$-vector (or a $\{0,\pm 1\}$-vector), $x$ of length $d+1$, such that the entries of $Px$ are the eigenvalues of $X$. If the entries of $Px$ are even, then $P x\equiv 0\pmod 2$, but $P$ is invertible modulo $2$ and therefore $x$ is zero modulo 2 implying that $X$ is empty.
\qed

We now turn to perfect state transfer on translation graphs. It will be useful to know the time at which perfect state transfer occurs in a graph. The following lemma follows directly from a theorem of Coutinho \cite{Coutinho2014}. Note that Coutinho's theorem holds for integer-weighted graphs, and so our lemma holds for integral signed translation graphs.

\begin{lem}[{\cite[Theorem 2.4.4]{Coutinho2014}}]\label{lem:timeofpst}
  Let $X$ be an integral translation graph with distinct eigenvalues $\theta_0>\cdots>\theta_n$ and define
  \[\delta:=\gcd\{\theta_0-\theta_r:r=1,\dots,n\}.\]
  Then, if $X$ admits perfect state transfer, it occurs at time $\pi/\delta$.\qed
\end{lem}

The next lemma uses the characters of $G$ to characterize perfect state transfer in translation graphs.

\begin{lem}
  \label{lem:psttranslation}
  Let $X=X(G,\C)$ be an integral translation graph and let $\delta$ be the greatest common divisor of its eigenvalue differences. Then, there is perfect state transfer from $0$ to $c$ in $X$ at time $\pi/\delta$, if and only if for each character, $\chi$ of $G$, we have
  \[\chi(c) = (-1)^{(|\C|-\chi(\C))/\delta}.\]
\end{lem}
\proof
  For each character, $\chi$ of $G$, define the $|G|\times |G|$ matrix $E_\chi$ by
  \[(E_\chi)_{g,h} := \frac1{|G|}\chi(h-g) = \frac1{|G|}\chi(h)\overline{\chi(g)},\]
  for all $g,h\in G$. The columns of $E_\chi$ are eigenvectors of $A:=A(X)$ with eigenvalue $\chi(\C)$ and we have
  \[A = \sum_\chi \chi(\C)E_\chi.\]
  This is a refinement of the spectral decomposition of $A$, and if $U(t)$ is the transition matrix at time $t$, then
  \[U(t) = \sum_\chi e^{it\chi(\C)}E_\chi.\]
  In particular,
  \begin{align*}
    U(t)_{0,c} & = \sum_\chi e^{it\chi(\C)}(E_\chi)_{0,c}\\
    & = \frac1{|G|}\sum_\chi e^{it\chi(\C)}\chi(c)\overline{\chi(0)}\\
    & = \frac1{|G|}\sum_\chi e^{it\chi(\C)}\chi(c).
  \end{align*}
  Since each term in the sum has absolute value one, we see that $|U(t)_{0,c}|=1$ if and only if the terms are all equal. This is equivalent to
  \[e^{it|\C|} = e^{it\chi(\C)}\chi(c)\]
  for all characters $\chi$ of $G$. By Lemma \ref{lem:timeofpst}, perfect state transfer must occur at time $\pi/\delta,$ and now the result follows.
\qed

It follows from Lemma \ref{lem:psttranslation} that perfect state transfer in a translation graph always occurs between $0$ and some element, say $a$, of order two. In this case, perfect state transfer also occurs between $g$ and $g+a$ for all $g\in G$.

We conclude this section by showing how, under certain conditions, we can partition the connection set of a translation graph and relate perfect state transfer on the original graph to perfect state transfer on Cayley graphs on the subsets.

\begin{lem}\label{lem:partition}
  Let $\C$ be an inverse-closed subset of the abelian group $G$ with a partition into inverse-closed subsets $\D_1$ and $\D_2$. Suppose $X(G,\D_2)$ is periodic with period $\tau$. Then $X(G,\D_1)$ admits perfect state transfer at time $\tau$ if and only if $X(G,\C)$ admits perfect state transfer at time $\tau$.
\end{lem}
\proof
  Let $A, A_1$ and $A_2$ be the adjacency matrices of $X=X(G,\C), X(G,\D_1)$ and $X(G,\D_2)$, respectively. Then $A= A_1+A_2$ and since $G$ is abelian, $A_1$ and $A_2$ commute. Therefore,
  \begin{align*}
    U_A(t) & = \exp(itA)\\
    & =\exp(it(A_1+A_2))\\
    & = \exp(itA_1)\exp(itA_2)\\
    & = U_{A_1}(t)U_{A_2}(t),
  \end{align*}
  for $t\in \R$.
  Since $X(G,\D_2)$ is periodic at time $\tau$, this implies $U_{A_2}(\tau)= I$, and so $U_A(\tau) = U_{A_1}(\tau)$. The result follows.
\qed

\section{Decomposing 2-circulants}

Recall that a $2$-circulant is a Cayley graph for an abelian group with a cyclic Sylow-$2$-subgroup. In this section we will see how we can decompose $2$-circulants into edge-disjoint subgraphs in a natural way. We can further decompose these subgraphs and we conclude the section with a theorem relating perfect state transfer on the big graph to perfect state transfer on the smaller graphs.

Let $G$ be an abelian group of order $2^dm$ where $m$ is odd, and assume $G$ has a cyclic Sylow-2-subgroup. Then $G\cong \Z_{2^d}\times H$, where $H\leq G$ is an abelian group of order $m$.
Let $X=X(G,\C)$ be an integral $2$-circulant. We partition $\C$ in the following way. Let $\C_k$ denote the subset of $\C$ consisting of elements that have order $2^km'$ for some odd number $m'$. It is clear that $\C_0,\dots,\C_d$ is a partition of $\C$. Further, $\C_k$ is power-closed (in particular inverse-closed) for all $k$.
We will show that the Cayley graph $X(G,\C_k)$ is a direct product of two Cayley graphs, one of which has complete bipartite graphs as its components. First, we need two lemmas.

\begin{lem}\label{lem:compbip}
  The Cayley graph $X(\Z_{2^d},[1])$ is isomorphic to the complete bipartite graph $K_{2^{d-1},2^{d-1}}$, for all $d\geq 1$.
\end{lem}
\proof
  We see that $[1]$ consists of all the odd numbers, so the even and odd numbers form a bipartition of $X(\Z_{2^d},[1])$ and the rest is clear.
\qed
Note that in a cyclic group, $[1]$ is the set of elements that generate the whole group. The next lemma is easily proved and this is left as an exercise for the reader.

\begin{lem}\label{lem:dirprod}
  Let $X = X(G,\C)$ be a Cayley graph. If $G \cong H_1\times H_2$ and $\C = \D_1\times \D_2$, then $X$ is the direct product of $Y_1=X(H_1,\D_1)$ and $Y_2=X(H_2,\D_2)$.\qed
\end{lem}

We can now prove the following theorem.
\begin{thm}\label{thm:dirprod}
  Let $X = X(G,\C)$ be an integral signed $2$-circulant. Suppose the order of $G$ is $2^dm$ where $m$ is odd and $d\geq 1$, and let $H$ be the unique subgroup of $G$ of order $m$. For $k\geq 0$, let $\C_k\subseteq \C$ be the set of elements in $\C$ of order $2^km'$ for some odd number $m'$.
  Further, let $K(d,k)$ denote the graph on $2^d$ vertices whose components are isomorphic to the complete bipartite graph $K_{2^{k-1},2^{k-1}}$, with the convention that $K(d,0)$ is the graph with adjacency matrix $I$. Then
  \[X(G,\C_k)\cong K(d,k)\times Y_k,\]
  where $Y_k=X(H,\D)$ is an integral signed Cayley graph of $H$. 
\end{thm}
\proof
  First, since $G\cong \Z_{2^d}\times H$, we can write every element of $G$ as $(x,h)$ with $x\in \Z_{2^d}$ and $h\in H$. Let $\C_k^+:=\C_k\cap\C^+$ and $\C_k^-:=\C_k\cap\C^-$. Since $X$ is an integral signed Cayley graph, $\C^+$ and $\C^-$ are power-closed, thus $\C_k^+$ and $\C_k^-$ are power-closed.
  We want to show that there are subsets $\D_1\seq \Z_{2^d}$ and $\D_2^+,\D_2^-\seq H$ such that $\C_k^+ = \D_1\times \D_2^+$ and $\C_k^-=\D_1\times \D_2^-$.
  We will show that if $(x_1,h_1),(x_2,h_2)\in \C_k$ (respectively $\C_k^+,\C_k^-$) then $(x_1,h_2),(x_2,h_1)\in \C_k$ (respectively $\C_k^+,\C_k^-$).

  Note that if $\langle (x_1,h_1)\rangle = \langle (x_2,h_2)\rangle\leq G$, then $\langle x_1\rangle = \langle x_2\rangle\leq \Z_{2^d}$ and $\langle h_1\rangle = \langle h_2\rangle\leq H.$
  Let $(x_1,h_1),(x_2,h_2)\in \C_k$. Then $x_1$ and $x_2$ have the same order, $2^k$, and since they are contained in the cyclic group $\Z_{2^d}$, we have $\langle x_1\rangle = \langle x_2\rangle$.
  But then we have $(x_1,h_2)\in [(x_2,h_2)]\seq\C_k$ and $(x_2,h_1)\in [(x_1,h_1)]\seq \C_k$, implying that $\C_k = \D_1\times \D_2$ for some $\D_1,\D_2$.

  Similarly we get $\C_k^+=\D_1'\times \D_2'$ and $\C_k^-=\D_1''\times \D_2''$. Since $\C_k$ is a disjoint union of $\C_k^+ $ and $\C_k^-$, and the elements of $\D_1$ all generate the same subgroup, we must have $\D_1'=\D_1''=\D_1$ and $\D_2'\cup\D_2''=\D_2$ with $\D_2'$ and $\D_2''$ disjoint.
  Thus we have found sets $\D_1\seq \Z_{2^d}$ and $\D_2^+,\D_2^-\seq H$ such that $\C_k^+ = \D_1\times \D_2^+$ and $\C_k^-=\D_1\times \D_2^-$.

  Then by Lemma \ref{lem:dirprod}, we have
  \begin{align*}
    X(G,\C_k)&\cong X(\Z_{2^d},\D_1)\times X(H,\D_2),\\
    X(G,\C_k^+)&\cong X(\Z_{2^d},\D_1)\times X(H,\D_2^+),\\
    X(G,\C_k^-)&\cong X(\Z_{2^d},\D_1)\times X(H,\D_2^-).
  \end{align*}
  Clearly, $\D_2^+$ and $\D_2^-$ are power-closed, so $Y_k:=X(H,\D_2)$ is an integral signed Cayley graph of $H$. The signed adjacency matrix of $X(G,\C_k)$ is
  \[A(X(\Z_{2^d},\D_1))\otimes (A(X(H,\D_2^+))-A(X(H,\D_2^-))).\]

  Further, we see that $\D_1$ consists of all the elements in $\Z_{2^d}$ of order $2^k$, implying that the components of $X(\Z_{2^d},\D_1)$ are isomorphic to $X(\Z_{2^k},[1])$, unless $k=0$, in which case it is the graph with a loop on each vertex and no other edges. Now the result follows from Lemma \ref{lem:compbip}.
\qed


\begin{cor}\label{cor:periodic}
  For $k\geq 3,$ the signed graphs $X(G,\C_k)$ are periodic at time $\pi/2.$
\end{cor}
\proof
  The eigenvalues of $X_k:=X(G,\C_k)$ are $\theta\lambda$ where $\theta$ and $\lambda$ are eigenvalues of $K(d,k)$ and the signed adjacency matrix of $Y_k$, respectively. The eigenvalues of $K(d,k)$ are $\pm 2^{k-1}$ and $0$ and so it follows that all the eigenvalues of $X_k$ are divisible by $2^{k-1}$. In particular, if $k\geq 3$ they are divisible by four. Then, looking at the spectral decomposition of the adjacency matrix of $X_k$, it is clear that $U(\pi/2) = I$, so $X_k$ is periodic at time $\pi/2$ for $k\geq 3$.
\qed

Let $X=X(G,\C)$ be an integral (signed) $2$-circulant where $G$ has order $2^dm$ with $m$ odd, and define $Y_k$ as in Theorem \ref{thm:dirprod} for $k=0,\dots,d$. Let $A_k$ denote the (signed) adjacency matrix of $Y_k$. We want to construct the adjacency matrix of $X$ using the matrices $A_0,\dots, A_d$.

Define the matrix $\M(A_0,\dots,A_d)$ recursively by
\[\M(A_0,A_1) :=\pmat{A_0&A_1\\A_1&A_0},\]
and for $k\geq 2$,
\begin{align*}
  \M(A_0,\dots,A_k) &:= I_2\otimes \M(A_0,\dots,A_{k-1})+\pmat{0&1\\1&0}\otimes(J_{2^{k-1}}\otimes A_k)\\
  & = \M(\M(A_0,\dots,A_{k-1}),J_{2^{k-1}}\otimes A_k).
\end{align*}
We claim that $\M(A_0,\dots,A_d)$ is the adjacency matrix of $X$. Recall that the adjacency matrix of the direct product of two graphs is the Kronecker product of their adjacency matrices. It is therefore easy to see that if $d=1$, then
\[A(X) = A(X(G,\C_0)) + A(X(G,\C_1)) = I_2\otimes A_0 + A(K_2)\otimes A_1 = \M(A_0,A_1).\]
(If $X$ is signed, $A(X)$ is the signed adjacency matrix of $X$.) Now the claim can be proved by induction.

If $X$ and $Y$ are graphs on the same vertex set with adjacency matrices $A$ and $B$, respectively, we define their \textit{semidirect product}, $X\ltimes Y$, as the graph with adjacency matrix $\M(A,B)$. This graph product is explored in a paper by Coutinho and Godsil \cite{Coutinho2018}, where they prove the following theorem.

\begin{thm}[{\cite[Theorem 5.2]{Coutinho2018}}]\label{thm:semidirect}
  Given graphs $X$ and $Y$ on the same vertex set $V$, with $A = A(X)$ and $B = A(Y)$, the graph $X\ltimes Y$ on vertex set $\{0,1\}\times V$ admits perfect state transfer if and only if one of the following holds.
  \begin{enumerate}[label=(\roman*)]
    \item For some $\tau\in\R^+, \lambda\in\C$ and $u\in V$, the matrices $A + B$ and $A-B$ are periodic at $u$ at time $\tau$ with respective phase factors $\lambda$ and $-\lambda$. In this case, perfect state transfer is between $(0,u)$ and $(1,u)$.
    \item For some $\tau\in\R^+, \lambda\in\C$ and $u\in V$, the matrices $A + B$ and $A-B$
    admit $uv$-perfect state transfer at time $\tau$ with the same phase factor $\lambda$. In this case, perfect state transfer is between $(0,u)$ and $(0,v)$, and between $(1,u)$ and $(1,v)$.
    \item For some $\tau\in\R^+, \lambda\in\C$ and $u\in V$, the matrices $A + B$ and $A-B$
    admit $uv$-perfect state transfer at time $\tau$ with respective phase factors $\lambda$ and $-\lambda$. In this case, perfect state transfer is between $(0, u)$ and $(1, v)$, and between $(1, u)$ and $(0, v)$.\qed
  \end{enumerate}
\end{thm}
Theorem \ref{thm:semidirect} is stated only for simple graphs in the original paper, but the proof also works if $X$ and $Y$ are signed graphs or if they have loops.

\section{Reducing to a simpler case}
In this section we will see how we can reduce the question of perfect state transfer on $2$-circulants to the case where the group has order $4m$ with $m$ odd.

Let $G$ be an abelian group of order $2^dm$ with $m$ odd, having a cyclic Sylow-2-subgroup. Note that if $d=0$, then $G$ has odd order and its Cayley graphs cannot have perfect state transfer by Theorem \ref{thm:pstvxtr}. Let us consider the case where $d=1$.

\begin{thm}\label{thm:mK2}
  Let $X=X(G,\C)$ be a loopless, signed Cayley graph for an abelian group $G$ of order $2m$, where $m$ is odd. If $X$ admits perfect state transfer, then $X \cong mK_2$ and the minimal time at which perfect state transfer occurs is $\pi/2$.
\end{thm}
\proof
  Suppose that (the signed adjacency matrix of) $X$ admits perfect state transfer at time $\tau$. Let $H$ be the unique subgroup of $G$ of order $m$, define $\C_0, \C_1, Y_0$ and $Y_1$ as in Theorem \ref{thm:dirprod}, and let $A_0,A_1$ be the (signed) adjacency matrices for $Y_0,Y_1$, respectively. Then $X\cong Y_0\ltimes Y_1$ and so we can apply Theorem \ref{thm:semidirect}.
  Recall that $Y_0$ and $Y_1$ are Cayley graphs for $H$, so $A_0+A_1$ and $A_0-A_1$ are weighted adjacency matrices of Cayley graphs for a group of odd order. Therefore, they do not admit perfect state transfer and so, (ii) and (iii) in Theorem \ref{thm:semidirect} cannot hold. So $A_0+A_1$ and $A_0-A_1$ are periodic at time $\tau$, with phase factors, say, $\lambda$ and $-\lambda$, respectively.
  Then $U_{A_0+A_1}(\tau) = \lambda I$ and $U_{A_0-A_1}(\tau) = -\lambda I$, thus we get
  \[U_{2A_1}(\tau) = U_{A_0+A_1}(\tau)U_{A_1-A_0}(\tau) = U_{A_0+A_1}(\tau)\overline{U_{A_0-A_1}(\tau)} = -I.\]
  This implies that if $\theta_0,\dots,\theta_n$ are the eigenvalues of $A_1$, then $e^{2i\tau\theta_k}=-1$ for all $k=0,\dots,n$. So, there are odd integers $m_k$ such that $2\tau\theta_k=m_k\pi$ and it follows that for all $j,k=0,\dots n$ we have $\theta_k/\theta_j =m_k/m_j$ with $m_k$ and $m_j$ both odd.

  Recall that a non-empty integral signed translation graph for a group of odd order has an odd eigenvalue, by Lemma \ref{lem:oddeval}. Furthermore, if such a graph is loopless, it must have even degree, which is then an even eigenvalue of the graph and it can be easily verified that this also holds in the signed case. Now, our graph $Y_1$ is either a (possibly signed) integral translation graph without loops, or its adjacency matrix may be written as $A_1 = A_1'\pm I$ where $A_1'$ is the (possibly signed) adjacency matrix of an integral translation graph without loops, and so in both cases, unless it is empty, $Y_1$ contains eigenvalues, $\theta, \theta'$ of different parities.
  But then, either the denominator or the enumerator of the reduced fraction $\theta/\theta'$ has to be even, so we conclude that $Y_1$ must be empty with loops, i.e., $A_1=\pm I$. The possible signing will not affect the rest of the proof, so we will drop it.

  We now have two possibilities: either $Y_0$ is empty in which case $X\cong mK_2$ as required, or $X$ is a Cartesian product, $K_2\square Y_0$ and $A(X) = A(K_2)\otimes I + I\otimes A_0$. Assume for contradiction the latter. Then
  \[U_{K_2\square Y_0}(\tau) =U_{K_2}(\tau)\otimes U_{Y_0}(\tau)= \gamma P,\]
  where $P$ is a permutation matrix with zero diagonal and $|\gamma|=1$.
  Then, both $U_{K_2}(\tau)$ and $U_{Y_0}(\tau)$ must be scalar multiples of permutation matrices, but since $Y_0$ is a Cayley graph for a group of odd order it does not have perfect state transfer and so we must have $U_{Y_0}(\tau)=\gamma'I$ for some scalar $\gamma'$.

  It is known (and easy to verify) that $K_2$ has perfect state transfer with minimal time $\pi/2$ and so $\tau = (2k+1)\pi/2$ for some integer $k$. Let $\theta_0',\dots,\theta_n'$ be the eigenvalues of $Y_0$. Then the eigenvalues of $U_{Y_0}(\tau)$ are $e^{i\pi(2k+1)\theta_j'/2}$.
  Since $Y_0$ has both even and odd eigenvalues, $U_{Y_0}(\tau)$ will have eigenvalues in both $\{\pm1\}$ and $\{\pm i\}$ and can therefore not be a scalar multiple of $I$.

  We have reached a contradiction and conclude that $X\cong mK_2$, having perfect state transfer with minimal time $\pi/2$.
\qed

Now recall that loops on a Cayley graph do not affect the existence of perfect state transfer, and so the only other Cayley graph of $G$ with perfect state transfer is the perfect matching with a loop on each vertex.
We have therefore covered the case where $d$ is at most one, and now we start considering the general case.

\begin{thm}\label{thm:timeofpst}
  Let $X=X(G,\C)$ be a 2-circulant. If perfect state transfer occurs on $X$, it occurs at time $\pi/2$.
\end{thm}
\proof
  Let $G$ have order $2^dm$ where $m$ is odd and suppose perfect state transfer occurs on $X$ at time $\tau$. We have seen what happens for $d=0,1$, so we assume that $d\geq 2.$

  Let $\C_0,\dots, \C_d$ be as before, define the graphs $Y_k$ as in Theorem \ref{thm:dirprod} and let $A_k$ be their adjacency matrices. Recall that
  \[A(X) = \M(A_0,\dots,A_d) = \M(\M(A_0,\dots,A_{d-1}),J_{2^{d-1}}\otimes A_d).\]
  We will use theorem \ref{thm:semidirect} on the matrices $\M(A_0,\dots,A_{d-1})$ and $J_{2^{d-1}}\otimes A_d$.  We know that perfect state transfer must occur between $0$ and the unique element of order two in $G$. This element is contained in the subgroup of $G$ generated by $\C_0\cup \C_1.$ Therefore, part (ii) of Theorem \ref{thm:semidirect} must apply and we have perfect state transfer at time $\tau$ on the matrix
  \[\M(A_0,\dots,A_{d-1}) - J_{2^{d-1}}\otimes A_d = \M(A_0 -A_d,\dots,A_{d-1}-A_d).\]
  Now, applying the theorem repeatedly, we get that perfect state transfer occurs on the matrix $\M(A_0-A_2,A_1-A_2)$ at time $\tau.$ This is a signed adjacency matrix of a Cayley graph for a group of order $2m$ and so by Theorem \ref{thm:mK2}, we have $\tau=\pi/2$ as required.
\qed

We are now ready for the main theorem of this section.

\begin{thm}\label{thm:main}
  Suppose $G$ is abelian of order $2^d m$ where $m$ is odd and $d\geq 2$,  and assume that the Sylow-$2$-subgroup of $G$ is cyclic. Let $G'$ denote the unique subgroup of $G$ with order $4m$. Then the integral Cayley graph $X(G,\C)$ admits perfect state transfer if and only if $X(G',G'\cap\C)$ admits perfect state transfer.
\end{thm}

\proof
  Note first that $X(G', G'\cap\C)$ admits perfect state transfer if and only if $X(G, G'\cap\C)$ admits perfect state transfer, since the components of $X(G,G'\cap\C)$ are all isomorphic to $X(G', G'\cap\C)$. Further, by Theorem \ref{thm:timeofpst}, if perfect state transfer occurs on any Cayley graph of $G$, it occurs at time $\pi/2$.

  Now define $\C_k$ for $k=0,\dots,d$ as before. We see that $G'\cap\C = \C_0\cup\C_1\cup\C_2$ and recall that by Corollary \ref{cor:periodic}, the graphs $X(G,\C_k)$ are periodic at time $\pi/2$, for $k\geq 3$. Then the result follows from Theorem \ref{lem:partition}.
\qed

Theorem \ref{thm:main} shows that the question of perfect state transfer on $2$-circulants can be reduced to groups that are isomorphic to $\Z_4\times H$, where $H$ has odd order. We will therefore devote Section \ref{sec:order4} to such groups.

Before we move on to the next section however, we have one more lemma which we will state in this section, since it applies to all abelian groups with a cyclic Sylow-2-subgroup of order at least four. Observe that any such group has a unique element of order two, and a unique pair of inverse elements of order four.

\begin{lem}\label{lem:parities}
  Let $G$ be an abelian group with a Sylow-2-subgroup that is cyclic and has order at least four. Let $\C$ be a power-closed subset such that $0\not\in \C$ and define $\C_k$ as before, for $k=0,1,2$. Let $a$ be the unique element of order two and $b,-b$ the unique pair of elements of order four. Then
  \begin{enumerate}[label = (\alph*)]
    \item $|\C_0|$ and $|\C_2|$ are even,
    \item $|\C_1|$ is odd if and only if $a\in \C$, and
    \item $|\C_2|$ is divisible by four if and only if $b\not\in \C$.
  \end{enumerate}
\end{lem}
\proof
  Since $\C_0,\C_1$ and $\C_2$ are power-closed, they are also inverse-closed. The only element in $G$ that is its own inverse is $a$, and if $a\in \C$, it is in $\C_1$. Therefore, the elements of $\C_0,\C_1,\C_2$ come in pairs, with the possible exception of $a$, from which (a) and (b) follow.

  For (c), write $G=\Z_{2^d}\times H$ where $H$ has odd order and observe that every element in $\C_2$ can be written of the form $c+h$ where $c\in \{b,-b\}$ and $h\in H$. Further, the elements $c+h,-c+h, c-h$ and $-c-h$ all generate the same subgroup, and provided that $h\neq 0$, they are all distinct. Therefore, with the exception of $b,-b$, we can partition $\C_2$ into subsets of size four in this way, and now (c) follows.
\qed

\section{Groups of order $4m$ with $m$ odd}\label{sec:order4}

In this section, $G$ is an abelian group of order $4m$ where $m$ is odd, having a cyclic Sylow-2-subgroup. We will give a characterization of connection sets $\C$ such that $X(G,\C)$ admits perfect state transfer. In the rest of the paper we will assume that $0\not\in \C$ (recall that this does not affect the existence of perfect state transfer).
We can write $G\cong \Z_4\times H$ where $H$ has order $m$. If $\C$ is a power-closed subset of $G$, we have $\C=\C_0\cup\C_1\cup\C_2,$ where $\C_k$ is defined as before. Let $a$ be the unique element of order two and $b,-b$ the unique pair of elements of order four.

We start with a simple necessary condition on $\C$.

\begin{lem}\label{lem:pstimplies2}
  Let $G = \Z_4\times H$ where $H$ has odd order $m$ and let $\C$ be an inverse-closed subset such that the Cayley graph $X=X(G,\C)$ admits perfect state transfer. Then either $a$ or $b$ is in $\C$, but not both.
\end{lem}
\proof
  Since $X$ has perfect state transfer, it is integral and so $\C$ is power-closed. Further, since $a$ is the unique element of order two in $G$, perfect state transfer must occur between $0$ and $a$ and by Theorem \ref{thm:timeofpst} it occurs at time $\pi/2$. Then, by Lemma \ref{lem:psttranslation} we have for each character $\chi$ of $G$,
  \begin{align}\label{eq:pstchar}
    \chi(a) = (-1)^{(|\C|-\chi(\C))/2}.
  \end{align}
  Let $\phi$ be the character of $G$ with $\phi(b) = i$ and $\phi(h) = 1$ for all $h\in H$. Then $\phi(a)=-1$, so by Equation (\ref{eq:pstchar}) we know that \((|\C|-\phi(\C))/2\)
  is odd. Further, we have $\phi(-b)=-i$, and now we see that
  \[\phi(\C_0) = |\C_0|,\quad \phi(\C_1) = -|\C_1|  \quad \text{and} \quad \phi(\C_2) = 0,\]
  and therefore
  \[|\C|-\phi(\C) = |\C|-\left(|\C_0|-|\C_1|\right) = 2|\C_1|+|\C_2|.\]
  Now, combining this with Lemma \ref{lem:parities}, we see that if $a$ and $b$ are both in $\C$, or if neither of them is in $\C$, then $2|\C_1|+|\C_2| = |\C|-\phi(\C)$ is divisible by four, contradicting that $(|\C|-\phi(\C))/2$ is odd. Therefore, we must have that $a\in\C$ or $b\in\C$, but not both.
\qed

Now, think of the elements of $G$ as pairs $(c,h)$ with $c\in \Z_4$ and $h\in H$. The set $\C_0$ consists only of elements $(0,h)$, and so we can view it as a subset of $H$. The set $\C_1$ has elements of the form $(a,h)$ for $h\in H$, and so there is a power-closed subset, $\C_1^*$ of $H$ such that $\C_1=\{a\}\times \C_1^*$, and moreover $0\in \C_1^*$ if and only if $a\in\C$.

Finally, the elements of $\C_2$ have the form $(\pm b,h)$ with $h\in H$ and since $\C_2$ is inverse-closed, $(b,h)\in \C_2$ if and only if $(-b,h)\in \C_2$.
Therefore there is a power-closed subset $\C_2^*$ of $H$ such that $\C_2=\{-b,b\}\times \C_2^*$ and $0\in\C_2^*$ if and only if $b\in\C$.

We can now prove another necessary condition for perfect state transfer to occur. We will later see that these two conditions along with integrality are also sufficient.

\begin{lem}\label{lem:pstimplies3}
  Let $G = \Z_4\times H$ where $H$ is an abelian group of odd order $m$ and let $\C$ be an inverse-closed subset such that the Cayley graph $X=X(G,\C)$ admits perfect state transfer. Then $\C_0=\C_1^*\diff\{0\} = \C_2^*\diff\{0\}$.
\end{lem}
\proof
  Since $X$ admits perfect state transfer, it must be integral, so $\C$ is power-closed. Further, by Lemma \ref{lem:pstimplies2}, either $a$ is in $\C$ or $b$ is in $\C$ but not both.

  Let $\phi_0,\phi_1,\phi_2$ be characters of $\Z_4$ with $\phi_0(b)=1$,  $\phi_1(b)=-1$ and $\phi_2(b)=i$ and let $\psi$ be an arbitrary character of $H$. Define the map $\chi_j:G\to \Cx^*$ by letting $\chi_j(c,h)=\phi_j(c)\psi(h)$ Then $\chi_j=\phi_j\psi$ is a character of $G$ for $j=0,1,2,$ and we have
  \begin{align*}
    \chi_j(\C) & =\chi_j(\C_0)+\chi_j(\C_1)+\chi_j(\C_2)\\
    & = \psi(\C_0) + \phi_j(a)\psi(\C_1^*) +(\phi_j(b)+\phi_j(-b))\psi(\C_2^*).
  \end{align*}
  Therefore
  \begin{align*}
    \phi_0\psi(\C) &= \psi(\C_0)+\psi(\C_1^*)+2\psi(\C_2^*),\\
    \phi_1\psi(\C) &= \psi(\C_0)+\psi(\C_1^*)-2\psi(\C_2^*),\\
    \phi_2\psi(\C) &= \psi(\C_0)-\psi(\C_1^*).
  \end{align*}

  Since $X$ has perfect state transfer, we have
  \begin{align*}
    \chi(a) &= (-1)^{(|\C|-\chi(\C))/2}
  \end{align*}
  for every character $\chi=\phi\psi$ of $G$.
  This means that whenever $\phi(b)=\pm i$, we have that $(|\C|-\chi(\C))/2$ is odd and otherwise, it is even. In other words, the following holds for all characters, $\psi$ of $H$:
  \begin{align}
    \frac{|\C|-(\psi(\C_0)+\psi(\C_1^*)+2\psi(\C_2^*))}2 & \quad\text{is even,}\label{first}\\[0.2cm]
    \frac{|\C|-(\psi(\C_0)+\psi(\C_1^*)-2\psi(\C_2^*))}2 & \quad\text{is even,}\label{second}\\[0.2cm]
    \frac{|\C|-(\psi(\C_0)-\psi(\C_1^*))}2 & \quad\text{is odd.}\label{third}
  \end{align}
  By adding (\ref{first}) and (\ref{second}) together, we find that
  $|\C|-\psi(\C_0)-\psi(\C_1^*)$ is even and by adding (\ref{second}) and (\ref{third}) we see that $|\C|-\psi(\C_0)+\psi(\C_2^*)$ is odd.

  If $a\not\in \C$, then $|\C|$ is even, $0\not\in\C_1^*$ and $0\in\C_2^*$, so this implies that $\psi(\C_0) +\psi(\C_1^*)$ and $\psi(\C_0)+\psi(\C_2^*\diff\{0\})$ are even for all characters, $\psi$ of $H$.
  If $a\in \C$, then we similarly get that $\psi(\C_0)+\psi(\C_1^*\diff\{0\})$ and $\psi(\C_0)+\psi(\C_2^*)$ are even for all $\psi$.

  Taking $\D:=\C_1^*\diff\{0\}$, we will consider the set $\S = (\C_0\diff\D)\cup(\D\diff\C_0)$. Our goal is to show that this set is empty, implying that $\C_0=\C_1^*\diff\{0\}$.

  Note that since $\C_0$ and $\D$ are both power-closed, their intersection is as well and then also their difference. Therefore $\S$ is a disjoint union of two power-closed sets and so it is itself power-closed. Further, we have
  \[\psi(\S)=\psi(\C_0\diff\D)+\psi(\D\diff\C_0) = \psi(\C_0)+\psi(\D)-2\psi(\C_0\cap\D)\]
  which by the above is even for all characters, $\psi$ of $H$.

  It follows that the Cayley graph $Y:=X(H,\S)$ has only even eigenvalues, but $Y$ is an integral Cayley graph for a group of odd order, so by Lemma \ref{lem:oddeval}, we conclude that $Y$ must be empty, implying that $\C_0=\C_1^*\diff\{0\}$.

  By taking $\D$ to be $\C_2^*\diff\{0\}$, we similarly get $\C_0=\C_2^*\diff\{0\}$, as desired.
\qed

We are now ready for the main theorem of this section.

\begin{thm}\label{thm:pstiff}
  Let $G\cong \Z_4\times H$ where $H$ is abelian of odd order. Then the Cayley graph $X(G,\C)$ admits perfect state transfer if and only if the following conditions hold:
  \begin{enumerate}[label = (\alph*)]
    \item $\C$ is power-closed,
    \item exactly one of $a$ and $b$ is in $\C$, and
    \item $\C_0=\C_1^*\diff\{0\} = \C_2^*\diff\{0\}$.
  \end{enumerate}
\end{thm}

\proof
  We have seen that a graph must be integral to have perfect state transfer, and that a translation graph is integral if and only if $\C$ is power-closed, thus perfect state transfer implies (a) and by lemmas \ref{lem:pstimplies2} and \ref{lem:pstimplies3}, it also implies (b) and (c).

  Suppose that conditions (a)-(c) of the theorem hold. By (a), $X$ is integral so by Lemma \ref{lem:psttranslation}, it suffices show that for every character $\chi$ of $G$ we have
  \begin{align}\label{eq:pstcharagain}
      \chi(a) = (-1)^{(|\C|-\chi(\C))/2}.
  \end{align}

  Note that every character $\chi$ of $G$ can be written uniquely as a product of two characters, $\chi=\phi\psi$ with $\phi$ a character of $\Z_4$ and $\psi$ a character of $H$, where $\chi(c,h)=\phi(c)\psi(h)$. Moreover, $\chi(a)=-1$ if and only if $\phi(a)=-1$, which holds if and only if $\phi(b)=\pm i$.

  Let $\chi=\phi\psi$ and assume first that $a\in \C$. Then $0\in\C_1^*$ and $0\not\in\C_2^*$, so by condition (c) we have $\C_0=\C_1^*\sm\{0\} = \C_2^*$. Therefore,
  \begin{align*}
    |\C| &= |\C_0|+|\C_1|+|\C_2|\\
    & = |\C_0| + |\C_1^*\sm\{0\}|+1+2|\C_2^*|\\
    & = 4|\C_0|+1.
  \end{align*}
  Further, we have
  \begin{align*}
    \chi(\C) &= \chi(\C_0)+\chi(\C_1)+\chi(\C_2)\\
    & = \psi(\C_0) + \phi(a)\psi\left(\C_1^*\sm\{0\}\right)+\phi(a) + (\phi(b)+\phi(-b))\psi(\C_2^*)\\
    & =\psi(\C_0) + \phi(a)\psi\left(\C_0\right)+\phi(a) + (\phi(b)+\phi(-b))\psi(\C_0)\\
    & = \psi(\C_0)\big(1+\phi(a)+\phi(b)+\phi(-b)\big) +\phi(a).\tag{$*$}
  \end{align*}
  Now if $\phi(b)=\pm i$, then $\phi(-b) = -\phi(b)$ and $\phi(a) = -1$. We see that in this case, $(*)$ becomes $\phi(a)=-1$ and so
  \begin{align*}
    \frac{|\C|-\chi(\C)}{2} & = \frac{4|\C_0|+1-(-1)}2\\
    & = 2|\C_0|+1
  \end{align*}
  which is odd and therefore, Equation (\ref{eq:pstcharagain}) holds.
  If $\phi(b) = -1$ then $\phi(-b)=-1$ and $\phi(a) = 1$ and so $(*)$ becomes $\phi(a)=1$ and
  \[\frac{|\C|-\chi(\C)}{2} = \frac{4|\C_0|}2 = 2|\C_0|,\]
  which is even, so again Equation (\ref{eq:pstcharagain}) holds in this case.
  Finally if $\phi(b) = 1$, then $(*) = 4\psi(\C_0)+\phi(a)$, so
  \[\frac{|\C|-\chi(\C)}{2} = \frac{4|\C_0|-4\psi(\C_0)}2 = 2\big(|\C_0|-\psi(\C_0)\big),\]
  which is even, as required.

  Then suppose $b\in \C$, and therefore $a\not\in \C$. Then condition (c) gives $\C_0=\C_1^* = \C_2^*\sm\{0\}$, and so
  \begin{align*}
    |\C| &= |\C_0|+|\C_1|+|\C_2|\\
    & = |\C_0| + |\C_1^*|+2|\C_2^*\sm\{0\}| +2\\
    & = 4|\C_0|+2.
  \end{align*}
  Moreover,
  \begin{align*}
    \chi(\C) &= \chi(\C_0)+\chi(\C_1)+\chi(\C_2)\\
    & = \psi(\C_0) + \phi(a)\psi\left(\C_1^*\right) + (\phi(b)+\phi(-b))\psi(\C_2^*\sm\{0\})+\phi(b)+\phi(-b)\\
    & =\psi(\C_0) + \phi(a)\psi\left(\C_0\right) + (\phi(b)+\phi(-b))\psi(\C_0)+\phi(b)+\phi(-b)\\
    & = \psi(\C_0)\big(1+\phi(a)+\phi(b)+\phi(-b)\big) +\phi(b)+\phi(-b).\tag{$**$}
  \end{align*}
  If $\phi(b)=\pm i$, then $(**)$ is zero, so
  \begin{align*}
    \frac{|\C|-\psi(\C)}2 &= \frac{4|\C_0|+2}2 = 2|\C_0|+1,
  \end{align*}
  which is odd. If $\phi(b)=-1$ then $(**) = -2$ and so
  \[\frac{|\C|-\psi(\C)}2 = \frac{4|\C_0|+2-2}2 = 2|\C_0|,\]
  which is even. Finally if $\phi(b) = 1$ then $(**) = 4\psi(\C_0)+2$ so
  \begin{align*}
    \frac{|\C|-\psi(\C)}2 &= \frac{4|\C_0|+2-(4\psi(\C_0)+2)}2\\
    & = \frac{4(|\C_0|-\psi(\C_0))}2\\
    & = 2(|\C_0|-\psi(\C_0)),
  \end{align*}
  again even. In each case, Equation (\ref{eq:pstcharagain}) holds, and we have shown that $X$ admits perfect state transfer.
\qed

We conclude this section with a fun corollary of Lemma \ref{lem:pstimplies3} about the spectrum of $X$.

\begin{cor}
  Let $G\cong \Z_4\times H$ where $H$ is abelian with odd order and suppose that the Cayley graph $X=X(G,\C)$ admits perfect state transfer. If $a\in \C$, then the eigenvalues of $X$ are all odd, and $-1$ is an eigenvalue with multiplicity $|G|/2$. If $a\not\in \C$, all eigenvalues of $X$ are even and $0$ is an eigenvalue with multiplicity $|G|/2$.
\end{cor}
\proof
  Since $X$ admits perfect state transfer, we have $\C_0=\C_1^*\diff\{0\}=\C_2^*\diff\{0\}$. Let $\phi_0,\phi_1,\phi_2,\phi_3$ be the four characters of $\Z_4$ defined uniquely by
  \[\phi_0(b)=1,\quad \phi_1(b)=-1,\quad \phi_2(b) =i,\quad \phi_3(b)=-i.\]
  From the proof of Lemma \ref{lem:pstimplies3}, the eigenvalues of $X$ are given by
  \begin{align*}
    \phi_0\psi(\C) &= \psi(\C_0)+\psi(\C_1^*)+2\psi(\C_2^*)=
    \begin{cases}
      4\psi(\C_0)+1 & \text{if }a\in\C\\
      4\psi(\C_0)+2 & \text{otherwise}
    \end{cases}\\[0.2cm]
    \phi_1\psi(\C) &= \psi(\C_0)+\psi(\C_1^*)-2\psi(\C_2^*) =
    \begin{cases}
      1 & \text{if }a\in\C\\
      -2 & \text{otherwise}
    \end{cases}\\[0.2cm]
    \phi_2\psi(\C)&=\phi_3\psi(\C) = \psi(\C_0)-\psi(\C_1^*) =
    \begin{cases}
      -1 & \text{if }a\in\C\\
      0 & \text{otherwise}
    \end{cases}
  \end{align*}
  for characters $\psi$ of $H$. We see that all eigenvalues are odd if $a\in \C$ and even otherwise. Further, since there are $|H|=|G|/4$ characters of $H$, and each of them gives two characters of $G$ with eigenvalue $0$ or $-1$, the result follows.
\qed

\section{Back to 2-circulants}

It now remains to combine Theorems \ref{thm:main} and \ref{thm:pstiff} to prove our characterization of connection sets that yield perfect state transfer. We will rephrase our conditions slightly here. Let $k\in \N$ and let $\S$ be a subset of an abelian group $G$ written additively. Then, we define the set $k\S:=\{k\cdot g:g\in \S\}$.

\begin{thm}\label{thm:final}
  Let $G$ be an abelian group of order $2^dm$ where $m$ is odd and suppose it has a cyclic Sylow-2-subgroup. Let $X=X(G,\C)$ be a Cayley graph. If $d=0$, there is no perfect state transfer on $X$ and if $d=1$, there is perfect state transfer if and only if $X$ is a matching. Suppose $d\geq 2$, let $a$ be the unique element of order two, and $b,-b$ the unique pair of elements of order four. Denote by $\C_k$ the set of elements in $\C$ with order $2^km'$ where $m'$ is odd. Then $X$ has perfect state transfer if and only if
  \begin{enumerate}[label = (\alph*)]
    \item $\C$ is power-closed,
    \item either $a$ or $b$ is in $\C$ but not both,
    \item $\C_0=4(\C_2\diff\{-b,b\})$, and
    \item $\C_1\diff\{a\} = 2(\C_2\diff\{-b,b\})$.
  \end{enumerate}
\end{thm}
\proof
  If $d=0$, then $G$ has odd order so $X$ has no perfect state transfer and for   $d=1$, the statement follows from Theorem \ref{thm:mK2}.

  Now suppose $d\geq 2$ and let $G'$ be the unique subgroup of $G$ of order $4m$. Then $G'\cong \Z_4\times H$ where $H$ has odd order, and by Theorem \ref{thm:main}, $X$ admits perfect state transfer if and only if $Y:=X(G', G'\cap\C)$ does. Note that $G'\cap\C=\C_0\cup\C_1\cup\C_2$.
  Further, $a,b\in G'$ and $a$ is the unique element of $G'$ of order two and $b,-b$ the unique pair of elements of $G'$ of order four.

  Define $\C_1^*,\C_2^*$ as before, so $\C_1=\{a\}\times\C_1^*$ and $\C_2=\{-b,b\}\times\C_2^*$. By Theorem \ref{thm:pstiff}, $Y$ has perfect state transfer if and only if conditions (a) and (b) hold and $\C_0=\C_1^*\diff\{0\}=\C_2^*\diff\{0\}.$

  Recall that $H$ has odd order. Then, for any $h\in H$, the elements $h$, $2h$ and $4h$ all generate the same cyclic subgroup because $2$ and $4$ are coprime to the order of $\langle h\rangle$. Thus, if $h\in \C_0$, we also have $2h,4h\in \C_0$ and same holds for $\C_1^*$ and $\C_2^*$.
  Therefore,
  \begin{align*}
    4(\C_2\diff\{-b,b\}) & = \{-4b,4b\}\times\{4h:h\in\C_2^*\diff\{0\}\}\\
    & = \{0\}\times(\C_2^*\diff\{0\}),
  \end{align*}
  and
  \begin{align*}
    2(\C_2\diff\{-b,b\}) & = \{-2b,2b\}\times\{2h:h\in\C_2^*\diff\{0\}\}\\
    & = \{a\}\times(\C_2^*\diff\{0\}).
  \end{align*}
  Now it is clear that condition (c) is equivalent to $\C_0=\C_2^*\diff\{0\}$ and condition (d) is equivalent to $\C_1^*\diff\{0\}=\C_2^*\diff\{0\}$ and the result follows.
\qed

\begin{xmpl}
Let $G=\Z_4\times\Z_3\times\Z_3$ with generators $b, g_1 $ and $g_2$ and $a=2b$. Let
\begin{align*}
  \C_0 & = \{g_1,\,2g_1,\,g_2,\,2g_2\}\\
  \C_1 & = \{a+ g_1,\,a+2g_1,\,a+g_2,\,a+2g_2\}\\
  \C_2 & = \{g_1\pm b,\,2g_1\pm b,\,g_2\pm b,\,2g_2\pm b\}
\end{align*}
and define $\C:=\C_0\cup\C_1\cup\C_2\cup\{a\}$ and $\C':=\C_0\cup\C_1\cup\C_2\cup\{-b,b\}$. We see that $4\C_2 = 2\C_1 = \C_0$ and so the graphs $X=X(G,\C)$ and $X'=X(G,\C')$ have perfect state transfer from $0$ to $a$. It can be shown by considering the automorphism group of these graphs that they are not circulants, and so we have found examples that do not arise from Ba\v si\'c's characterization.

The graph $X$ is shown in Figure \ref{fig:pstgraph} with the vertices $0$ and $a$ shown in white.
\begin{figure}[h!]
  \centering
  \includegraphics[scale=0.8]{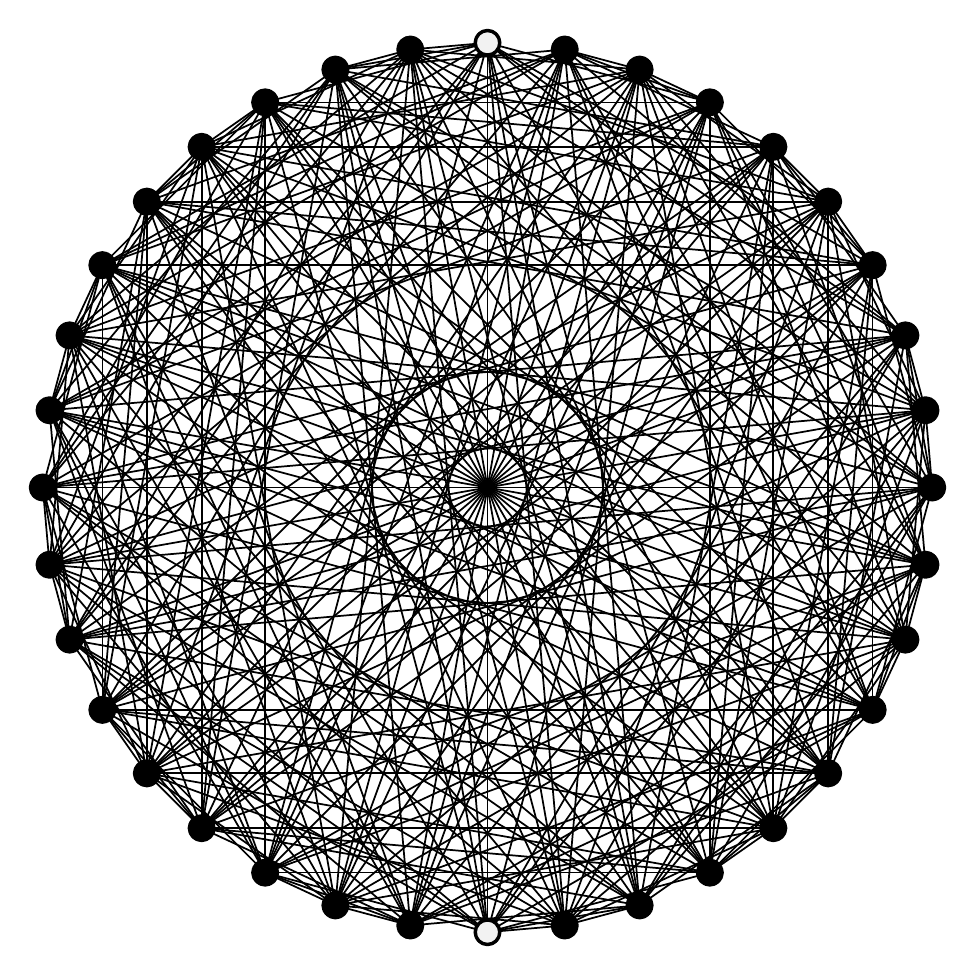}
  \caption{Cayley graph of $\Z_4\times\Z_3\times\Z_3$ with perfect state transfer}
  \label{fig:pstgraph}
\end{figure}
\end{xmpl}

\section{Further work}
In this paper, we have considered Cayley graphs for groups of the form $G=\Z_{2^d}\times H$ where $H$ is an arbitrary abelian group of odd order, and our main result is a generalization of Ba\v si\'c's result where $G$ (equivalently $H$) is cyclic.

A general question that we might ask is when do we get perfect state transfer on vertex transitive graphs? What about Cayley graphs? The problem here is that we use the fact that the characters of an abelian group are the eigenvectors of its Cayley graphs, and from this we get a characterization of perfect state transfer on translation graphs.

This does not hold in general for an arbitrary Cayley graph of a non-abelian group, but similar things do happen in \textit{normal Cayley graphs}, that is Cayley graphs for which the connection set is conjugacy-closed. But then we run into other problems. Our proofs rely heavily on the fact that if perfect state transfer occurs in our graphs, the minimal time at which it occurs is $\pi/2$. This is not the case for all Cayley graphs, in fact, on Cayley graphs of elementary abelian 2-groups, perfect state transfer can happen arbitrarily fast \cite{chan2013complex}. Such groups are certainly worth exploring, but it is unlikely that the techniques we have used in this paper will be fruitful for this.

Then suppose $G$ is a (not necessarily abelian) group of order $2^dm$ where $m$ is odd and suppose it has a cyclic Sylow-$2$-subgroup, $S$. It can be shown that in this case, $G$ has a unique subgroup $H$ of order $m$, and so this subgroup is normal in $G$. It follows that $G$ is a semidirect product of $S$ and $H$. Now we ask two questions.

\begin{itemize}
  \item If perfect state transfer occurs on a Cayley graph of a group that has a cyclic Sylow-2-subgroup, does it necessarily occur at time $\pi/2$ at the earliest?
  \item Can we characterize perfect state transfer in normal Cayley graphs of such groups?
\end{itemize}

\vspace{0.6cm}

\bibliographystyle{plain}



\end{document}